\documentclass{jctart08}
\usepackage{graphics, graphicx}     

\usepackage{cite}
\usepackage{amsmath}
\usepackage{amssymb}
\usepackage{amsthm}
\theoremstyle{plain}
\newtheorem{Theorem}{Theorem}[section] %
\newtheorem{Lemma}{Lemma}[section]

\theoremstyle{definition}

\renewcommand{\thetable}{\thesection.\arabic{table}}

\addto\captionsrussian{\def\refname{References}}
\newenvironment{Proof} 
{\par\noindent{\it Proof of}} 
{\hfill$\vspace{5mm}\scriptstyle\blacksquare$} 

\numberwithin{equation}{section} 
\numberwithin{figure}{section} 
\numberwithin{table}{section} 

\begin{document}

\setcounter{page}{1}

\markboth{M.I. Isaev, R.G. Novikov}{Stability estimates for determination of potential from the impedance boundary map}

\title{Stability estimates for determination of potential from the impedance boundary map}
\date{}
\author{ { M.I. Isaev and R.G. Novikov}}

\maketitle
{\bf Abstract}
\begin{abstract}
We study the impedance boundary map (or Robin-to-Robin map) for the Schr\"odinger equation
in open bounded domain at fixed energy in multidimensions. We give global stability estimates for
determining potential from these boundary data and, as corollary, from the Cauchy data set. Our results
include also, in particular, an extension of the Alessandrini identity to the case of the impedance boundary map.
\end{abstract}

\section{Introduction}

We consider the  Schr\"odinger equation
\begin{equation}\label{eq} 
	-\Delta \psi  + v(x)\psi = E\psi, \ \  x \in  D, \, E\in\mathbb{R},
\end{equation}
where
\begin{equation}\label{eq_c}
	\begin{aligned}
		 D &\text{ is an open bounded domain in } \mathbb{R}^d,\ d\geq 2, \ \\ &\text{with } \partial D \in C^2,
 	\end{aligned}
\end{equation}
\begin{equation}\label{eq_c2}
 v \in \mathbb{L}^{\infty}(D).
\end{equation}
We consider the impedance boundary map $\hat{M}_\alpha = \hat{M}_{\alpha,v}(E)$ defined by
\begin{equation}\label{def_M}
	\hat{M}_\alpha [\psi]_\alpha = [\psi]_{\alpha-\pi/2}
\end{equation}
 for all sufficiently regular solutions $\psi$ of  equation (\ref{eq}) in $\bar{D} = D \cup \partial D$, where
\begin{equation}\label{[psi]_alpha}
	[\psi]_\alpha = [\psi(x)]_\alpha = \cos\alpha  \, \psi(x)  -  \sin\alpha\, \frac{\partial \psi}{\partial \nu}|_{\partial D}(x),\ \  
	x \in \partial D, \ \alpha \in \mathbb{R}
\end{equation}
and $\nu$ is the outward normal to $\partial D$. One can show(see Lemma \ref{Number_alpha}) that there is not more than a countable number of $\alpha\in\mathbb{R}$ such that $E$ is an eigenvalue for the operator $-\Delta + v$ in $D$
with the boundary condition 
\begin{equation}\label{bound_cond}
 \cos\alpha  \, \psi|_{\partial D}  -  \sin\alpha\, \frac{\partial \psi}{\partial \nu}|_{\partial D}  = 0.
\end{equation} 
Therefore, for any energy level $E$ we can assume that  for some fixed  $\alpha\in\mathbb{R}$
\begin{equation}\label{correct_M}
	\begin{aligned}
		\text{ $E$ is not an eigenvalue for the operator $-\Delta + v$ in $D$} \\
		\text{with boundary condition (\ref{bound_cond})}
	\end{aligned}
\end{equation}
and, as a corollary, $\hat{M}_\alpha$ can be defined correctly.

Note that the impedance boundary map $\hat{M}_\alpha$ is reduced to the 
Dirichlet-to-Neumann(DtN) map if $\alpha=0$ and is reduced to the 
Neumann-to-Dirichlet(NtD) map if $\alpha=\pi/2$. 
The map $\hat{M}_\alpha$ can be called also as the Robin-to-Robin map.
General Robin-to-Robin map was considered, in particular, in \cite{GM2009}.

We consider the following inverse boundary value problem for equation (\ref{eq}).

{\bf Problem 1.1.}
 Given $\hat{M}_\alpha$ for some fixed $E$ and $\alpha$, find $v$. 

This problem can be considered as the Gel’fand inverse boundary value problem for the Schr\"odinger equation at fixed energy
(see \cite{Gelfand1954}, \cite{Novikov1988}). At zero energy this problem can be considered also as a generalization of the Calderon problem of the electrical impedance tomography (see \cite{Calderon1980}, \cite{Novikov1988}).

Problem 1.1 includes, in particular, the following questions: (a) uniqueness, (b) reconstruction,
(c) stability.

Global uniqueness theorems and global reconstruction methods for Problem 1.1 with $\alpha =0$
were given for the first time in \cite{Novikov1988} in dimension $d \geq 3$  and in \cite{Buckhgeim2008} in dimension $d=2$.

Global stability estimates for Problem 1.1 with $\alpha =0$ were given for the first time in \cite{Alessandrini1988} 
in dimension $d \geq 3$ and in \cite{NS2010} in dimension $d=2$.  
A principal improvement of the result of \cite{Alessandrini1988} was given recently in \cite{Novikov2011} (for the zero energy case).
Due to 
  \cite{Mandache2001} these logarithmic 
stability results are optimal (up to the value of the exponent). 
An extention of the instability estimates of  \cite{Mandache2001}  to the case of the non-zero energy as well as to the case 
 of Dirichlet-to-Neumann map given on the energy intervals was given in \cite{IsaevDtN}.

Note also that for the Calderon problem (of the electrical impedance tomography) in its initial formulation the global uniqueness 
was firstly proved in \cite{SU1987} for $d\geq 3$ and in \cite{Nachman1996} for $d=2$.
 
 It should be noted that in  most of previous works on inverse boundary value problems for equation (\ref{eq}) at fixed $E$
 it was assumed in one way or another that 
 $E$ is not a Dirichlet eigenvalue for the operator $-\Delta + v$ in $D$, see \cite{Alessandrini1988}, \cite{Mandache2001}, \cite{Novikov1988}, \cite{Novikov2011}, \cite{NS2010}, \cite{NS2011}, \cite{S2011}. Nevertheless, the results of \cite{Buckhgeim2008} can be
 considered as global uniqueness  and reconstruction results for Problem 1.1 in dimension $d=2$ with general $\alpha$.
 
 In the present work we give global stability estimates for Problem 1.1 in dimension $d\geq 2$ with general $\alpha$. 
 These results are presented in detail in Section 2. 
 
 In addition, in the present work we establish some basic properties of the impedance boundary map with general $\alpha$. In particular, we extend 
 the Alessandrini identity to this general case. These results are presented in detail in Section 3.
	
	In a subsequent paper we plan to give also global reconstruction method for Problem 1.1 in multidimensions with general $\alpha$.
\section{Stability estimates}
In this section we always assume that $D$ satisfies (\ref{eq_c}). 

We will use 
the fact that  if $v_1$, $v_2$ are potentials satisfying (\ref{eq_c2}), (\ref{correct_M})
for some fixed $E$ and $\alpha$, then 
\begin{equation}
	\hat{M}_{\alpha,v_1}(E) - \hat{M}_{\alpha,v_2}(E) \text { is a bounded operator in } \mathbb{L}^{\infty}(\partial D),
\end{equation}
where $\hat{M}_{\alpha,v_1}(E)$, $\hat{M}_{\alpha,v_2}(E)$ denote the impedance boundary maps for $v_1$, $v_2$, respectively.
Actually, under our assumptions, $\hat{M}_{\alpha,v_1}(E) - \hat{M}_{\alpha,v_2}(E)$ is a compact operator in $\mathbb{L}^{\infty}(\partial D)$
(see Corollary 3.1). 

Let 
\begin{equation}\label{norm_A}
	\begin{aligned}
		||A|| \text{ denote the norm of an operator}\\
		A: \mathbb{L}^{\infty}(\partial D) \rightarrow \mathbb{L}^{\infty}(\partial D).
	\end{aligned}
\end{equation}

Let the Cauchy data set ${\cal C}_v$  for equation (\ref{eq}) be defined by:
\begin{equation}\label{eq_v}
	{\cal C}_v = \left\{\left(\psi|_{\partial D}, \frac{\partial \psi}{\partial \nu}|_{\partial D}\right) :
	\begin{array}{r}
		\text{ for all sufficiently regular solutions $\psi$ of }\\
		\text{ equation (\ref{eq}) in $\bar{D} = D \cup \partial D$}
	\end{array}
	\right\}.
\end{equation}
In addition, the Cauchy data set ${\cal C}_v$ can be represented as the graph of the impedance boundary map $\hat{M}_\alpha = \hat{M}_{\alpha,v}(E)$
defined by (\ref{def_M}) under assumptions (\ref{correct_M}).

\subsection{Estimates for $d\geq 3$}

In this subsection we assume for simplicity that
\begin{equation}\label{assumption}
		 v \in W^{m,1}(\mathbb{R}^d) \text{ for some } m > d, \ \mbox{supp}\, v \subset D,\ 
\end{equation}
where
\begin{equation}\label{2.5}
	W^{m,1}(\mathbb{R}^d) = \{v:\  \partial^J v \in L^1(\mathbb{R}^d),\  |J| \leq m \},\ m\in \mathbb{N}\cup 0,
\end{equation}
where 
\begin{equation}
J \in (\mathbb{N}\cup 0)^d,\ |J| = \sum\limits_{i=1}\limits^{d}J_i,\ \partial^J v(x) 
= \frac{\partial^{|J|} v(x)}{\partial x_1^{J_1}\ldots \partial x_d^{J_d}}.
\end{equation}
Let
\begin{equation}
	||v||_{m,1} = \max\limits_{|J|\leq m} ||\partial^J v||_{L^1(\mathbb{R}^d)}.
\end{equation}
Note also that (\ref{assumption}) $\Rightarrow$ (\ref{eq_c2}).
\begin{Theorem}\label{Theorem_2.1} Let $D$ satisfy (\ref{eq_c}), where $d \geq 3$. 
Let  $v_1$, $v_2$ satisfy  (\ref{assumption}) and (\ref{correct_M}) 
for some fixed $E$ and $\alpha$.
Let $||v_{j}||_{m,1} \leq N,\  j = 1,2$, for some $N > 0$. 
Let $\hat{M}_{\alpha,v_1}(E)$ and $\hat{M}_{\alpha,v_2}(E)$ denote 
the impedance boundary maps for
$v_1$ and  $v_2$, respectively. Then
\begin{equation}\label{estimation_A}
	||v_1 - v_2||_{L^\infty(D)} \leq 
	C_\alpha \left(\ln\left(3+\delta_\alpha ^{-1}\right)\right)^{-s}, \ \ 0<s\leq (m-d)/m,
\end{equation}
where $C_\alpha = C_\alpha(N,D,m,s,E)$,  $\delta_\alpha  = ||\hat{M}_{\alpha,v_1}(E) - \hat{M}_{\alpha,v_2}(E)||$ 
is defined according  to (\ref{norm_A}).
\end{Theorem}

\noindent
{\bf Remark 2.1.} Estimate (\ref{estimation_A})
with $\alpha = 0$ is a variation of the result of \cite{Alessandrini1988} (see also \cite{Novikov2011}). 

Proof of Theorem \ref{Theorem_2.1} is given in Section 5. This proof is based on results presented in Sections 3, 4.

Theorem \ref{Theorem_2.1} implies the following corollary:

\noindent
{\bf Corollary 2.1.} { \it Let $D$ satisfy (\ref{eq_c}), where $d \geq 3$. Let potentials $v_1$, $v_2$ satisfy  (\ref{assumption}). Then
\begin{equation}\label{estimation_cor}
	||v_1 - v_2||_{L^\infty(D)} \leq 
	\min\limits_{\alpha \in \mathbb{R}} 
	C_\alpha \left(\ln\left(3+\delta_\alpha^{-1}\right)\right)^{-s}, \ \ 0<s\leq (m-d)/m,
\end{equation} 
where $C_\alpha$  and $\delta_\alpha$ at fixed $\alpha$ are the same that in Theorem \ref{Theorem_2.1}.
}

Actually, Corollary 2.1 can be considered as global stability estimate for determining potential $v$ from its Cauchy data set ${\cal C}_v$ for 
equation (\ref{eq}) at fixed energy $E$, where $d\geq 3$.

\subsection{Estimates for $d = 2$}

In this subsection we assume for simplicity that
\begin{equation}\label{assumption2}
		 v \in C^{2}(\bar{D}),  \ \ \mbox{supp}\, v \subset D.
\end{equation}
Note also that (\ref{assumption2}) $\Rightarrow$ (\ref{eq_c2}).

\begin{Theorem}\label{Theorem_2.2}
Let $D$ satisfy (\ref{eq_c}), where $d = 2$.
Let  $v_1$, $v_2$ satisfy  (\ref{assumption2}) and (\ref{correct_M}) 
for some fixed $E$ and $\alpha$.
Let $||v_{j}||_{C^{2}(\bar{D})} \leq N,\  j = 1,2$, for some $N > 0$. 
Let $\hat{M}_{\alpha,v_1}(E)$ and $\hat{M}_{\alpha,v_2}(E)$ denote 
the impedance boundary maps for
$v_1$ and  $v_2$, respectively. Then
\begin{equation}\label{estimation_2}
	||v_1 - v_2||_{L^\infty(D)} \leq 
	C_\alpha \left(\ln\left(3+\delta_\alpha^{-1}\right)\right)^{-s}
	\left( \ln \left(3\ln\left(3+\delta_\alpha^{-1}\right)\right) \right)^2
	, \ \ 0<s\leq 3/4,
\end{equation}
where $C_\alpha = C_\alpha(N,D,s,E)$,   $\delta_\alpha  = ||\hat{M}_{\alpha,v_1}(E) - \hat{M}_{\alpha,v_2}(E)||$  is defined according  to (\ref{norm_A}).
\end{Theorem}

\noindent
{\bf Remark 2.2.} Theorem \ref{Theorem_2.2} for $\alpha = 0$ was given in \cite{NS2010} with $s = 1/2$ and in \cite{S2011}
 with $s = 3/4$.

Proof of Theorem \ref{Theorem_2.2} is given in Section 7. This proof is based on results presented in Sections 3, 6.

Theorem \ref{Theorem_2.2} implies the following corollary:

\noindent
{\bf Corollary 2.2.} { \it Let $D$ satisfy (\ref{eq_c}), where $d = 2$. Let potentials $v_1$, $v_2$ satisfy  (\ref{assumption2}). Then
\begin{equation}\label{estimation_cor2}
	||v_1 - v_2||_{L^\infty(D)} \leq 
	\min\limits_{\alpha \in \mathbb{R}} 
	C_\alpha \left(\ln\left(3+\delta_\alpha^{-1}\right)\right)^{-s}
	\left( \ln \left(3\ln\left(3+\delta_\alpha^{-1}\right)\right) \right)^2
	, \ \ 0<s\leq 3/4,
\end{equation} 
where $C_\alpha$ and $\delta_\alpha$ at fixed $\alpha$ are the same that in Theorem \ref{Theorem_2.2}.
}

Actually, Corollary 2.2 can be considered as global stability estimate for determining potential $v$ from its Cauchy data set ${\cal C}_v$ for 
equation (\ref{eq}) at fixed energy $E$, where $d =  2$.  

\subsection{Concluding remarks}
\ \ \ \ 
Theorems 2.1, 2.2 and Lemma \ref{Number_alpha} imply the following corollary:

\noindent
{\bf Corollary 2.3.} { \it  Under assumptions (1.2), (1.3), real-valued potential $v$
is uniquely determined by its Cauchy data ${\cal C}_v$ at fixed real energy $E$ }.

To our knowledge the result of Corollary 2.3 for $d\geq 3$ was not yet completely proved in the literature.

Let $\sigma_{\alpha,v}$ denote the spectrum of the operator $-\Delta+v$ in $D$ with boundary condition (\ref{bound_cond}). 

\noindent
{\bf Remark 2.3.} In Theorems \ref{Theorem_2.1} and \ref{Theorem_2.2} we do not assume that 
$E \notin \sigma_{\alpha,v_1}\cup \sigma_{\alpha,v_2}$ namely for $\alpha = 0$ in contrast with \cite{Alessandrini1988},
\cite{Novikov2011}, \cite{NS2010}, \cite{NS2011}, \cite{S2011}.
In addition, in fact, in Corollaries 2.1 and  2.2 there are no special assumptions on $E$ and $\alpha$ at all.
Actually, the stability estimates of  \cite{Alessandrini1988},
\cite{Novikov2011}, \cite{NS2010}, \cite{NS2011}, \cite{S2011} 
make no sense for $E \in \sigma_{0,v_1}\cup \sigma_{0,v_2}$ and  are too weak  if 
$dist (E, \sigma_{0,v_1}\cup \sigma_{0,v_2})$ is too small.

\noindent
{\bf Remark 2.4.}		 
The stability estimates of Subsections 2.1 and 2.2 admit principal improvement in the sense
described in \cite{Novikov2011}, \cite{NN2009}, \cite{S2011N}. In particular, Theorem \ref{Theorem_2.1} with $s=m-d$ (for $d=3$ and $E=0$)
follows from results presented in Sections 3, 4 of the present work and results presented in Section 8 of \cite{Novikov2011}.
In addition, estimates (2.8), (2.9) for $s =(m-d)/d$ admit a proof technically very similar to the proof
of Theorem 2.1, presented in Section 5. Possibility of such a proof of estimate (2.8) for $s=(m-d)/d$, $\alpha=0$, $E=0$
was mentioned, in particular, in \cite{Palamodov2011}.

\noindent
{\bf Remark 2.5.}
The stability estimates of Subsections 2.1 and 2.2 can be extended  to the case when we do not assume that 
		 $\mbox{supp}\,v\subset D$ or, by other words, that $v$ is zero near the bounadry.
		 In this connection see, for example, \cite{Alessandrini1988}, \cite{NS2010}.
		 
In the present work we do not develop Remarks 2.4 and 2.5 in detail because of restrictions in time.

Note also that Theorems \ref{Theorem_2.1} and \ref{Theorem_2.2} remain valid with complex-valued potentials $v_1$, $v_2$ 
and complex $E$, $\alpha$. Finally, we note that in Theorems 2.1, 2.2 and Corollaries 2.1, 2.2 with real $\alpha$, constant  $C_\alpha$
can be considered as independent of $\alpha$.



\section{Some basic properties of the impedance boundary map}

\begin{Lemma}\label{Lemma_3.1}
	Let $D$ satisfy (\ref{eq_c}). Let potential $v$ satisfy (\ref{eq_c2}) and (\ref{correct_M}) 
	for some fixed $E$ and $\alpha$. Let $\hat{M}_\alpha = \hat{M}_{\alpha,v}(E)$ denote 
the impedance boundary map for $v$. Then
\begin{equation}\label{M+I}
	\begin{aligned}
		&\left(\sin\alpha \,\hat{M}_\alpha + \cos\alpha \,\hat{I}\right) [\psi]_\alpha = \psi|_{\partial D},\\
		&\left(\cos\alpha \, \hat{M}_\alpha - \sin\alpha \,\hat{I}\right) [\psi]_\alpha = \frac{\partial\psi}{\partial\nu}|_{\partial D},
	\end{aligned}
\end{equation} 

\begin{equation}\label{Sym_M}
	\int\limits_{\partial D}
	[\psi^{(1)}]_\alpha  \hat{M}_\alpha [\psi^{(2)}]_\alpha  d x= 
		\int\limits_{\partial D}
	[\psi^{(2)}]_\alpha \hat{M}_\alpha [\psi^{(1)}]_\alpha d x
\end{equation}
 for all sufficiently regular solutions $\psi$, $\psi^{(1)}$, $\psi^{(2)}$ of  equation (\ref{eq}) in $\bar{D}$, 
 where $\hat{I}$ denotes the identity operator on $\partial D$ and $[\psi]_\alpha$ is defined by (\ref{[psi]_alpha}).
\end{Lemma}
Note that identities (\ref{M+I}) imply that
\begin{equation}\label{M_alpha_beta}
	\left(\sin(\alpha_1-\alpha_2)\hat{M}_{\alpha_1} + \cos(\alpha_1-\alpha_2)\hat{I}\right) 
	\left(\sin(\alpha_2-\alpha_1)\hat{M}_{\alpha_2} + \cos(\alpha_2-\alpha_1)\hat{I}\right) = \hat{I},
\end{equation}
under the assumptions of Lemma \ref{Lemma_3.1} fulfilled simultaneously for $\alpha = \alpha_1$ and $\alpha = \alpha_2$.

Note also that from (\ref{Sym_M}) we have that
\begin{equation}\label{Sym_M++}
	\int\limits_{\partial D}
	[\phi^{(1)}]_\alpha  \hat{M}_\alpha [\phi^{(2)}]_\alpha  d x= 
		\int\limits_{\partial D}
	[\phi^{(2)}]_\alpha \hat{M}_\alpha [\phi^{(1)}]_\alpha d x
\end{equation}
for all sufficiently regular functions $\phi^{(1)}$, $\phi^{(2)}$ on $\partial D$.

\begin{Proof} {\it Lemma {\ref{Lemma_3.1}}}.
Identities (\ref{M+I}) follow from definition (\ref{def_M}) of the map  $\hat{M}_\alpha$.

To prove (\ref{Sym_M}) we use, in particular, the Green formula
\begin{equation}\label{G_F}
	\int\limits_{\partial D} \left( \phi^{(1)} \frac{\partial\phi^{(2)}}{\partial \nu} 
		- \phi^{(2)} \frac{\partial\phi^{(1)}}{\partial \nu} \right)  dx = 
		\int\limits_{D} \left(\phi^{(1)} \Delta \phi^{(2)} - \phi^{(2)} \Delta \phi^{(1)}\right) dx,
\end{equation}
where $\phi^{(1)}$ and $\phi^{(2)}$ are arbitrary sufficiently regular functions in $\bar{D}$.
Using (\ref{G_F}) and the identities
	\begin{equation}\label{identities}
		\psi^{(1)} \Delta \psi^{(2)}   
		= (v-E)\psi^{(1)}\psi^{(2)} = 
		\psi^{(2)} \Delta \psi^{(1)}  \ \ \text{ in $D$,}
	\end{equation}
we obtain that
	\begin{equation}\label{psi_1+2}
	\int\limits_{\partial D} \left( \psi^{(1)} \frac{\partial\psi^{(2)}}{\partial \nu} 
		- \psi^{(2)} \frac{\partial\psi^{(1)}}{\partial \nu} \right)  dx =  0.
	\end{equation}
	Using (\ref{psi_1+2}), we get that
	\begin{equation}\label{3.5}
		\begin{aligned}
		\int\limits_{\partial D} 
		\left( \cos\alpha\,\psi^{(1)} - \sin\alpha\frac{\partial\psi^{(1)}}{\partial \nu} \right)
		\left( \sin\alpha\,\psi^{(2)} + \cos\alpha\frac{\partial\psi^{(2)}}{\partial \nu} \right) dx = \\
		=\int\limits_{\partial D} 
			\left( \cos\alpha\,\psi^{(2)} - \sin\alpha\frac{\partial\psi^{(2)}}{\partial \nu} \right)
			\left( \sin\alpha\,\psi^{(1)} + \cos\alpha\frac{\partial\psi^{(1)}}{\partial \nu} \right) dx.
		\end{aligned}
	\end{equation}
	Identity (\ref{Sym_M}) follows from (\ref{3.5}) and definition (\ref{def_M}) of the map  $\hat{M}_\alpha$.
\end{Proof}

\begin{Theorem}\label{Th_ident}
	Let $D$ satisfy (\ref{eq_c}). Let two potentials $v_1$,  $v_2$   satisfy  (\ref{eq_c2}), (\ref{correct_M}) for 
	some fixed $E$ and  $\alpha$. 
	Let  $\hat{M}_{\alpha,v_1} = \hat{M}_{\alpha,v_1}(E)$, $\hat{M}_{\alpha,v_2} =\hat{M}_{\alpha,v_2}(E)$ 
	denote the impedance boundary maps 
	for  $v_1$, $v_2$, respectively. Then
	\begin{equation}\label{Ident}
			\int\limits_{D} \left(v_1 - v_2\right) \psi_1 \psi_2\, dx 
			=	\int\limits_{\partial D}
			[\psi_1]_\alpha\left(\hat{M}_{\alpha,v_1} - \hat{M}_{\alpha,v_2}\right) [\psi_2]_\alpha dx
	\end{equation}
	for all sufficiently regular  solutions $\psi_1$ and $\psi_2$ of equation (\ref{eq}) in  $\bar{D}$ with 
	$v = v_1$ and  $v= v_2$, respectively, where $[\psi]_\alpha$ is defined by (\ref{[psi]_alpha}).
\end{Theorem}

\begin{Proof} {\it Theorem {\ref{Th_ident}}}.
	As in (\ref{identities}) we have that
		\begin{equation}\label{3.7}
		\begin{aligned}
			\psi_1 \Delta \psi_2   
				= (v_2-E)\psi_1\psi_2, \\
		\psi_2 \Delta \psi_1  
				= (v_1-E)\psi_1\psi_2.
		\end{aligned}
	\end{equation}
	Combining (\ref{3.7}) with (\ref{G_F}),  (\ref{M+I}) and (\ref{Sym_M++}), we obtain that
	\begin{equation}
		\begin{aligned}
					\int\limits_{D} \left(v_1(x) - v_2(x)\right) \psi_1(x) \psi_2(x) dx =
				\int\limits_{\partial D} \left( \psi_2 \frac{\partial\psi_1}{\partial \nu} 
		- \psi_1 \frac{\partial\psi_2}{\partial \nu} \right)  dx = \\
		= \int\limits_{\partial D} 
		\left( \sin\alpha \, \hat{M}_{\alpha, v_2} + \cos \alpha \, \hat{I} \right) [\psi_2]_\alpha
		\left( \cos\alpha \, \hat{M}_{\alpha, v_1} - \sin \alpha \, \hat{I} \right) [\psi_1]_\alpha dx \ -\\
		- \int\limits_{\partial D} 
		\left( \sin\alpha \, \hat{M}_{\alpha, v_1} + \cos \alpha \, \hat{I} \right) [\psi_1]_\alpha
		\left( \cos\alpha \, \hat{M}_{\alpha, v_2} - \sin \alpha \, \hat{I} \right) [\psi_2]_\alpha  dx = \\
		= \int\limits_{\partial D}
			[\psi_1]_\alpha\left(\hat{M}_{\alpha, v_1} - \hat{M}_{\alpha, v_2}\right) [\psi_2]_\alpha  dx.
		\end{aligned}
	\end{equation}
\end{Proof}

\noindent
{\bf Remark 3.1.}
 Identity (\ref{Ident}) for  $\alpha = 0$  is reduced to Alessandrini's identity (Lemma 1 of \cite{Alessandrini1988}).

Let 
	$G_\alpha(x,y,E)$ be the Green function for the operator $\Delta - v + E$ in $D$
	 with the impedance boundary condition ({\ref{bound_cond}}) under assumptions (\ref{eq_c}),  (\ref{eq_c2}) and (\ref{correct_M}).
Note that  
\begin{equation}\label{3.10}
	G_\alpha(x,y,E) = G_\alpha(y,x,E),  \ \ \ x,y \in \bar{D}.
\end{equation}
The symmetry (\ref{3.10}) is proved in Section 9.

\begin{Theorem}\label{G+M}
	Let $D$ satisfy (\ref{eq_c}). Let potential $v$ satisfy (\ref{eq_c2}) and (\ref{correct_M}) 
	for some fixed $E$ and $\alpha$ such that $\sin\alpha \neq 0$. Let 
$G_\alpha(x,y,E)$ be the Green function for the operator $\Delta - v + E$ in $D$
	 with the impedance boundary condition ({\ref{bound_cond}}). 
	Then for $x,y\in \partial D$ 
	\begin{equation}\label{eq_G+M}
		M_\alpha(x,y,E) =  \frac{1}{\sin^2\alpha}\ G_\alpha(x,y,E) - \ctg \alpha\ \delta_{\partial D}(x-y),
	\end{equation}
	where  $M_\alpha(x,y,E)$ and  $\ \delta_{\partial D}(x-y)$ denote 
the Schwartz kernels of the impedance boundary map $\hat{M}_\alpha = \hat{M}_{\alpha,v}(E)$
 and the identity operator $\hat{I}$ on $\partial D$, respectively, 
 where $\hat{M}_\alpha$ and $\hat{I}$ are considered as linear integral operators.
\end{Theorem}
\begin{Proof} {\it Theorem \ref{G+M}.} 
	Note that
	\begin{equation}\label{3.15}
		[\phi]_{\alpha-\pi/2} = \frac{1}{\sin^2\alpha}\ \sin\alpha \,\phi|_{\partial D} - \ctg \alpha\ [\phi]_\alpha.
	\end{equation}
	for all suffuciently regular functions  $\phi$ in some neighbourhood of $\partial D$ in $D$.
	Since $G_\alpha$ is the Green function for equation (\ref{eq}) we have that
	\begin{equation}\label{3.12}
		\psi(y) = 	\int\limits_{\partial D} \left( \psi(x) \frac{\partial G_\alpha}{\partial \nu_x}(x,y,E) 
		-  G_\alpha(x,y,E)\frac{\partial\psi}{\partial \nu}(x) \right) dx, \ \  y\in D,
	\end{equation}
	for all  suffuciently regular  solutions $\psi$ of equation ({\ref{eq}}).
	Using (\ref{3.12}) and impedance boundary condition  ({\ref{bound_cond}}) for $G_\alpha$, we get that
	\begin{equation}\label{3.14}
		\begin{aligned}
		\sin\alpha \, \psi(y)  = 
		\sin\alpha \,
		\int\limits_{\partial D} \left( \psi(x) \frac{\partial G_\alpha}{\partial \nu_x}(x,y,E) 
		-  G_\alpha(x,y,E)\frac{\partial\psi}{\partial \nu}(x) \right) dx =\\
		=\int\limits_{\partial D} [\psi(x)]_\alpha   G_\alpha(x,y,E) dx,\ \  y\in D.
		\end{aligned}
	\end{equation}
	Due to (\ref{Sym_M++}) we have that
	\begin{equation}\label{Sym_M2}
		M_\alpha(x,y,E) = M_\alpha(y,x,E), \ \  x,y \in \partial D.
	\end{equation}
	Combining (\ref{def_M}), (\ref{3.15}), (\ref{3.14}) and (\ref{Sym_M2}),  we obtain (\ref{eq_G+M}).
\end{Proof}

\noindent
{\bf Corollary 3.1.} {\it Let assumtions of Theorem \ref{Th_ident} hold. Then
\begin{equation}\label{M_compact}
	\text{$\hat{M}_{\alpha,v_1}(E) - \hat{M}_{\alpha,v_2}$(E) is a compact operator in $\mathbb{L}^{\infty}(\partial D)$}.
\end{equation}}

\noindent
{\it Scheme of the proof of Corollary 3.1.} 
	Let 	$G_{\alpha, v_1 }(x,y,E)$ and $G_{\alpha, v_2}(x,y,E)$ 
	be the Green functions for the operator $\Delta  - v + E$ in $D$ with the impedance boundary condition ({\ref{bound_cond}}) 
	for $v = v_1$ and 	$v = v_2$, respectively. Using (\ref{3.10}), we find that
	\begin{equation}\label{3.17}
		\begin{aligned}
			G_{\alpha,v_1}(x,y,E) = \int\limits_{D} 
				 G_{\alpha,v_1}(x,\xi,E) \left( \Delta_\xi - v_2(\xi) + E\right) G_{\alpha,v_2}(\xi,y,E)\, d\xi,	\\
			G_{\alpha,v_2}(x,y,E) = \int\limits_{D} 
				\left( \Delta_\xi - v_1(\xi) + E\right) G_{\alpha,v_1}(x,\xi,E) G_{\alpha,v_2}(\xi,y,E)\, d\xi, \\
		\int\limits_{\partial D} \left( G_{\alpha,v_1}(x,\xi,E)  \frac{\partial G_{\alpha,v_2}}{\partial \nu_\xi} (\xi,y,E)
		- G_{\alpha,v_2}(\xi,y,E) \frac{\partial G_{\alpha,v_1}}{\partial \nu_\xi}(x,\xi,E)  \right)  d\xi = 0,\\
		x, y \in D. 
		\end{aligned}
	\end{equation}
	
	Combining (\ref{3.17}) with (\ref{G_F}), we get that
	\begin{equation}\label{3.19}
		\begin{aligned}
			G_{\alpha,v_1} (x,y,E) -  G_{\alpha,v_2}(x,y,E) =
							\int\limits_{D}
				\left( v_1(\xi) -v_2(\xi) \right) G_{\alpha,v_1}(x,\xi,E) G_{\alpha,v_2}(\xi,y,E)\, d\xi,\\
					x, y \in D. 
		\end{aligned}
	\end{equation}
	
	The proof of (\ref{M_compact}) for the case of $\sin \alpha \neq 0$  can be completed proceeding from
	(\ref{M_alpha_beta}), (\ref{eq_G+M}), (\ref{3.19}) and estimates of \cite{Lanzani2004} and \cite{Begehr2009} on
	$G_\alpha(x,y,E)$ for $v\equiv 0$.
	
	 Corollary 3.1 for the Dirichlet-to-Neumann case ($\sin\alpha = 0$)  was given in \cite{Novikov1988}.   
\hfill$\vspace{5mm}\scriptstyle\blacksquare$

\begin{Lemma}\label{Number_alpha}
Let $D$ satisfy (\ref{eq_c}).  Let  $v$ be a real-valued potential satisfying (\ref{eq_c2}). Then
for any fixed  $E\in \mathbb{R}$ there is not more than countable number of $\alpha \in \mathbb{R}$ such that $E$ is an eigenvalue for 
the operator $-\Delta+v$ in $D$ with boundary condition (\ref{bound_cond}).
\end{Lemma}
\begin{Proof} { \it Lemma \ref{Number_alpha}.} 
Let $\psi^{(1)}$, $\psi^{(2)}$ be eigenfunctions for 
the operator $-\Delta+v$ in $D$ with boundary condition (\ref{bound_cond}) 
for $\alpha = \alpha^{(1)}$ and $\alpha = \alpha^{(2)}$, respectively. Then	
	\begin{equation}
		\sin\left(\alpha^{(1)}  - \alpha^{(2)}\right) \int\limits_{\partial D} \psi^{(1)}\psi^{(2)} dx=
		\sin\alpha^{(1)} \sin\alpha^{(2)} 	
		\int\limits_{\partial D} \left( \psi^{(1)} \frac{\partial\psi^{(2)}}{\partial \nu} 
		- \psi^{(2)} \frac{\partial\psi^{(1)}}{\partial \nu} \right) dx = 0.
	\end{equation}
Since in the separable space $\mathbb{L}^2(\partial D)$ there is not more than countable 
orthogonal system of functions, we obtain the assertion of  Lemma \ref{Number_alpha}.
\end{Proof}

\noindent
{\bf Remark 3.1}
The assertion of Lemma \ref{Number_alpha} remains valid for the case of $\alpha \in \mathbb{C}$.

\section{Faddeev functions}
We consider the Faddeev functions $G$, $\psi$, $h$ (see \cite{Faddeev1965}, \cite{Faddeev1974}, \cite{Henkin1987}, \cite{Novikov 1988}):
\begin{equation}\label{4.1}
	\psi(x,k) = e^{ikx} + \int\limits_{\mathbb{R}^d} G(x-y,k) v(y)\psi(y,k) dy, 
\end{equation}
\begin{equation}\label{4.2}
	G(x,k) = e^{ikx} g(x,k), \ \ \ g(x,k) = - (2\pi)^{-d} \int\limits_{\mathbb{R}^d} \frac{e^{i\xi x} d\xi}{\xi^2 + 2k\xi},
\end{equation}
where $x\in \mathbb{R}^d$, $k\in \mathbb{C}^d$, $\mbox{Im}\, k\neq 0$, $d\geq 3$,
\begin{equation}\label{4.4}
	h(k,l) = (2\pi)^{-d} \int\limits_{\mathbb{R}^d} e^{-ilx}v(x) \psi(x,k) dx,
\end{equation}
where 
\begin{equation}\label{4.5}
	 k,l \in \mathbb{C}^d, \ k^2=l^2,\ \mbox{Im}\, k = \mbox{Im}\, l \neq 0.
\end{equation}
One can consider (\ref{4.1}), (\ref{4.4}) assuming that 
\begin{equation}\label{4.6}
	v \text{ is a sufficiently regular function on } \mathbb{R}^d
	\text{ with suffucient decay at infinity.}
\end{equation}
For example, in connection with Problem 1.1, one can consider (\ref{4.1}), (\ref{4.4}) assuming that
\begin{equation}\label{4.7}
	v \in \mathbb{L}^{\infty}(D), \ \ \ v \equiv 0 \text{ on } \mathbb{R}\setminus D.
\end{equation}

We recall that (see \cite{Faddeev1965}, \cite{Faddeev1974}, \cite{Henkin1987}, \cite{Novikov 1988}): 
\begin{itemize}
\item The function $G$ satisfies the equation
\begin{equation}\label{4.8}
	(\Delta+k^2) G(x,k) = \delta(x), \ \ x\in\mathbb{R}^d, \ \ k \in \mathbb{C}^d\setminus \mathbb{R}^d;
\end{equation}

\item Formula (\ref{4.1}) at fixed $k$ is considered as an equation for 
\begin{equation}
	\psi = e^{ikx}\mu(x,k),
\end{equation}
where $\mu$ is sought in $\mathbb{L}^{\infty}(\mathbb{R}^d)$; 

\item As a corollary of (\ref{4.1}), (\ref{4.2}), (\ref{4.8}), 
$\psi$ satisfies (\ref{eq}) for $E=k^2$; 

\item The Faddeev functions $G$, $\psi$, $h$ are (non-analytic) continuation to the complex domain of functions    
 of the classical scattering theory for the Schr\"odinger equation (in particular, $h$ is a generalized "`scattering"' amplitude). 
\end{itemize} 
 
 In addition, $G$, $\psi$, $h$ in their zero energy restriction, that is for $E=0$, were considered for the first time in \cite{Beals1985}.
The Faddeev functions $G$, $\psi$, $h$ were, actually, rediscovered in \cite{Beals1985}.

Let 
\begin{equation}
\begin{aligned}
		\Sigma_E = \left\{ k\in \mathbb{C}^d: k^2 = k_1^2 + \ldots + k_d ^2 = E\right\},\\
		\Theta_E = \left\{ k\in \Sigma_E,\  l\in\Sigma_E: \mbox{Im}\, k = \mbox{Im}\, l\right\}.
\end{aligned}		
\end{equation}
Under the assumptions of Theorem \ref{Theorem_2.1}, we have that:
\begin{equation}\label{mu_1}
	\mu(x,k) \rightarrow 1 \ \  \text{ as } \ \ |\mbox{Im}\,k|\rightarrow \infty 
\end{equation}
and, for any $\sigma>1$,
\begin{equation}\label{mu_2}
	|\mu(x,k)| + |\nabla \mu(x,k)|  \leq \sigma \ \  \text{ for } \ \ |\mbox{Im}\,k| \geq r_1(N,D,E,m,\sigma),
\end{equation}
where $x\in \mathbb{R}^d$, $k \in \Sigma_E$;
\begin{equation}\label{lim_1}
			\hat{v}(p) = \lim\limits_
		{\scriptsize
			\begin{array}{c}
			(k,l)\in \Theta_E,\, k-l=p\\
			|\mbox{Im}\,k|=|\mbox{Im}\,l|\rightarrow \infty
			\end{array}
		} h(k,l)\ \ \ 	 \text{ for any } p\in \mathbb{R}^d,
	\end{equation}
\begin{equation}\label{lim_2}
	\begin{aligned}
		|\hat{v}(p) - h(k,l)|\leq \frac{c_1(D,E,m)N^2}{\rho}  	\ \ \text{ for }  (k,l) \in \Theta_E, \ \  p = k-l,\\
		|\mbox{Im}\,k| = |\mbox{Im} \,l| = \rho \geq r_2(N,D,E,m), \\ p^2 \leq 4(E+\rho^2),
	\end{aligned}
\end{equation}
		where 
		\begin{equation}
			\hat{v}(p) = (2\pi)^{-d}\int\limits_{\mathbb{R}^d} e^{ipx} v(x)dx, \ \ p\in \mathbb{R}^d.
		\end{equation}

Results of the type (\ref{mu_1}) go back to \cite{Beals1985}.  Results of the type (\ref{lim_1}), (\ref{lim_2}) 
(with less precise right-hand side in (\ref{lim_2})) go back to \cite{Henkin1987}. In the present work estimate (\ref{mu_2}) 
is given according to \cite{Novikov1999}, \cite{Novikov2006}. Estimate (\ref{lim_2}) follows, for example, from the estimate
\begin{equation}\label{4.15.new}
	\begin{aligned}
	\| \Lambda^{-s} g(k) \Lambda^{-s}\|_{\mathbb{L}^2(\mathbb{R}^d)\rightarrow \mathbb{L}^2(\mathbb{R}^d)} = O(|k|^{-1}) \ \ \text{ as } \ |k|\rightarrow \infty,\\
		k\in \mathbb{C}^d\setminus \mathbb{R}^d,\ |k| = (|\mbox{Re}\,k|^2 +|\mbox{Im}\,k|^2)^{1/2},
	\end{aligned}
\end{equation}
for $s>1/2$, where $g(k)$ denotes the integral operator with the Schwartz kernel $g(x-y,k)$ and $\Lambda$ denotes the multiplication operator by the function $(1+|x|^2)^{1/2}$. Estimate (\ref{4.15.new}) was formulated, first, in \cite{LN1987} for $d\geq 3$. Concerning proof of (\ref{4.15.new}), see \cite{Weder1991}. 

In addition, we have that:
\begin{equation}\label{4.15}
	\begin{aligned}
	h_2(k,l) - h_1(k,l) = (2\pi)^{-d} \int\limits_{\mathbb{R}^d} \psi_1(x,-l) (v_2(x) - v_1(x)) \psi_2(x,k) dx 	\\ 	 
	\text{ for }  (k,l) \in \Theta_E,\ |\mbox{Im}\,k| = |\mbox{Im} \,l| \neq 0,\\ \text{ and $v_1$, $v_2$ satisfying (\ref{4.6}),}
	\end{aligned}
\end{equation}
\begin{equation}\label{4.16}
	\begin{aligned}
	h_2(k,l) - h_1(k,l) = (2\pi)^{-d} \int\limits_{\partial D} [\psi_1(\cdot,-l)]_\alpha 
	\left(\hat{M}_{\alpha,v_2} - \hat{M}_{\alpha,v_1}\right) 
	[\psi_2(\cdot,k)]_\alpha dx\\
	\text{ for }  (k,l) \in \Theta_E,\ |\mbox{Im}\,k| = |\mbox{Im} \,l| \neq 0,\\ \text{ and $v_1$, $v_2$ satisfying (\ref{correct_M}), (\ref{4.7}),}
	\end{aligned}
\end{equation}
where $h_j$, $\psi_j$ denote $h$ and $\psi$ of (\ref{4.4}) and (\ref{4.1}) for $v = v_j$, and $\hat{M}_{\alpha, v_j}$ denotes
the impedance boundary map of (\ref{def_M}) for $v = v_j$, where $j=1,2$.

Formula (\ref{4.15}) was given in \cite{Novikov1996}. Formula (\ref{4.16}) follows from Theorem \ref{Th_ident} and (\ref{4.15}).
Formula (\ref{4.16}) for $\alpha = 0$ was given in \cite{Novikov2005}.

\section{Proof of Theorem \ref{Theorem_2.1}}

Let 
\begin{equation}\label{6.3}
\begin{aligned}
	\mathbb{L}^{\infty}_{\mu}(\mathbb{R}^d) = \{u\in \mathbb{L}^{\infty}(\mathbb{R}^d): \|u\|_\mu<+\infty\},\\
	\|u\|_\mu = \mbox{ess}\,\sup\limits_{p\in\mathbb{R}^d} (1+|p|)^{\mu}|u(p)|, \ \ \ \mu>0.
\end{aligned}
\end{equation}

Note that 
\begin{equation}\label{6.4}
	\begin{aligned}
	w \in \mathbb{W}^{m,1}(\mathbb{R}^d) \Longrightarrow \hat{w} \in 	\mathbb{L}^{\infty}_{\mu}(\mathbb{R}^d)\cap {\cal C}(\mathbb{R}^d),\\
		\|\hat{w}\|_\mu \leq c_2(m,d) \|w\|_{m,1} \ \ \ \text{ for } \ \ \mu = m, 
	\end{aligned}
\end{equation}
where $\mathbb{W}^{m,1}$,  $\mathbb{L}^{\infty}_{\mu}$ are the spaces of (\ref{2.5}), (\ref{6.3}),
\begin{equation}
	\hat{w}(p) = (2\pi)^{-d}\int\limits_{\mathbb{R}^d} e^{ipx}w(x) dx, \ \ \ p\in \mathbb{R}^d.
\end{equation}

Using the inverse Fourier transform formula
\begin{equation}
	w(x) = \int\limits_{\mathbb{R}^d} e^{-ipx}\hat{w}(p) dp, \ \ \ x\in \mathbb{R}^d,
\end{equation}
we have that
\begin{equation}\label{5.4}
	\begin{aligned}
		\|v_1 - v_2\|_{\mathbb{L}^{\infty}(D)} \leq 
		\sup\limits_{x\in \bar{D}}|\int\limits_{\mathbb{R}^d} e^{-ipx}\left(\hat{v}_2(p) - \hat{v}_1(p)\right) dp| \leq\\
		\leq I_1(r) + I_2(r) \ \ \ \text{ for any }  \ r>0,
	\end{aligned}
\end{equation}
where
\begin{equation}\label{6.7}
	\begin{aligned}
		I_1(r) = \int\limits_{|p|\leq r} |\hat{v}_2(p) - \hat{v}_1(p)| dp, \\
		I_2(r) = \int\limits_{|p|\geq r} |\hat{v}_2(p) - \hat{v}_1(p)| dp.
	\end{aligned}
\end{equation}

Using (\ref{6.4}), we obtain that
\begin{equation}\label{6.8}
	|\hat{v}_2(p) - \hat{v}_1(p)| \leq 2c_2(m,d) N (1+|p|)^{-m}, \ \ \ p\in \mathbb{R}^d.
\end{equation}

Due to (\ref{lim_2}), we have that
\begin{equation}\label{6.9}
\begin{aligned}
	|\hat{v}_2(p) - \hat{v}_1(p)| \leq |h_2(k,l) - h_1(k,l)| + \frac{2c_1(D,E,m)N^2}{\rho},\\
	p\in \mathbb{R}^d,\   p = k-l,\ (k,l)\in \Theta_E,\\ |\mbox{Im}\,k| = |\mbox{Im}\,l| = \rho \geq r_2(N,D,E,m),\\
		p^2 \leq 4(E+\rho^2).
\end{aligned}
\end{equation}

Let 
\begin{equation}
	\begin{aligned}
	c_3 = (2\pi)^{-d}\int\limits_{\partial D} dx, \ \ \ L = \max\limits_{x\in \partial D}|x|,\\
	\delta_\alpha = \|\hat{M}_{\alpha, v_2}(E) - \hat{M}_{\alpha, v_1}(E)\|,
	\end{aligned}
\end{equation}
where $\|\hat{M}_{\alpha, v_2}(E) - \hat{M}_{\alpha, v_1}(E)\|$ is defined according to (\ref{norm_A}).

Due to (\ref{4.15}), (\ref{4.16}), we have that
\begin{equation}\label{6.11}
	\begin{aligned}
	|h_2(k,l) - h_1(k,l)| \leq  c_3
	\|[\psi_1(\cdot,-l)]_\alpha\|_{\mathbb{L}^{\infty}(\partial D)}\, \delta_\alpha \,
	\|[\psi_2(\cdot,k)]_\alpha\|_{\mathbb{L}^{\infty}(\partial D)},\\
	(k,l)\in \Theta_E, \ |\mbox{Im}\,k| = |\mbox{Im}\,l| \neq 0.
	\end{aligned}
\end{equation}
Using (\ref{[psi]_alpha}), (\ref{mu_2}), we find that
\begin{equation}\label{6.12}
\begin{aligned}
	\|[\psi(\cdot,k)]_\alpha\|_{\mathbb{L}^{\infty}(\partial D)} \leq
	c_4(E)\, \sigma \,\exp\bigg(|\mbox{Im}\,k|(L+1)\bigg), \\
	k \in \Sigma_E, \ |\mbox{Im}\,k|\geq r_1(N,D,E,m,\sigma).
\end{aligned}
\end{equation}
Here and bellow in this section the constant $\sigma$ is the same that in (\ref{mu_2}).
 
Combining (\ref{6.11}) and (\ref{6.12}), we obtain that
\begin{equation}\label{6.13}
	\begin{aligned}
		|h_2(k,l) - h_1(k,l)| \leq c_3 \left(c_4(E)\sigma\right)^2 \,\exp\bigg(2\rho(L+1)\bigg) \delta_\alpha, \\
		(k,l)\in\Theta_E, \ \rho = |\mbox{Im}\,k|= |\mbox{Im}\,l|\geq r_1(N,D,E,m,\sigma).
	\end{aligned}
\end{equation}
Using (\ref{6.9}), (\ref{6.13}), we get that
\begin{equation}\label{6.14}
	\begin{aligned}
		|\hat{v}_2(p) - \hat{v}_1(p)| \leq c_3 \left(c_4(E)\sigma\right)^2 \,\exp\bigg(2\rho(L+1)\bigg) \delta_\alpha 
		+ \frac{2c_1(D,E,m)N^2}{\rho},\\
		p\in\mathbb{R}^d,\  p^2 \leq 4(E+\rho^2),\  \rho \geq r_3(N,D,E,m,\sigma),
	\end{aligned}
\end{equation}
where $r_3(N,D,E,m,\sigma)$ is such that
\begin{equation}\label{5.14new.}	
	\rho \geq r_3(N,D,E,m,\sigma) \Longrightarrow 
	\left\{
	\begin{aligned}
	&\rho \geq r_1(N,D,E,m,\sigma), \\ &\rho \geq r_2(N,D,E,m), \\ &\rho^{2/m} \leq 4(E+\rho^2).
	\end{aligned}\right.
\end{equation}
Let 
\begin{equation}
	c_5 = \int\limits_{p\in \mathbb{R}^d,  |p|\leq 1} d p, \ \ \ c_6 =  \int\limits_{p\in \mathbb{R}^d,  |p|= 1} d p.
\end{equation}
Using (\ref{6.7}), (\ref{6.14}), we get that
\begin{equation}\label{5.15}
\begin{aligned}
	I_1(r) \leq c_5 r^d \left(c_3 \left(c_4(E)\sigma\right)^2 \,\exp\bigg(2\rho(L+1)\bigg) \delta_\alpha 
		+ \frac{2c_1(D,E,m)N^2}{\rho}\right),\\  r>0, \ r^2\leq 4(\rho^2 + E), \ \rho \geq r_3(N,D,E,m,\sigma).
\end{aligned}
\end{equation}
Using (\ref{6.7}), (\ref{6.8}), we find that for any $r>0$
\begin{equation}\label{5.16}
	\begin{aligned}
		I_2(r) \leq 2 c_2(m,d) N c_6 \int\limits_{r}\limits^{+\infty} \frac{dt}{t^{m-d+1}} \leq \frac{2c_2(m,D)Nc_6}{m-d}\frac{1}{r^{m-d}}.
	\end{aligned}
\end{equation}
 Combining  (\ref{5.4}), (\ref{5.15}), (\ref{5.16}) for $r = \rho^{1/m}$ and (\ref{5.14new.}), we get that 
 \begin{equation}\label{5.18}
 \begin{aligned}
 	\|v_1 - v_2\|_{\mathbb{L}^{\infty}(D)} \leq c_7(D,\sigma) \rho^{d/m} e^{2\rho(L+1)} \delta_\alpha + 
 	c_8(N,D,E,m)\rho^{-\frac{m-d}{m}}, \\
 	\rho \geq r_3(N,D,E,m,\sigma).
 	\end{aligned}
 \end{equation}
 
We fix some $\tau \in (0,1)$ and let
\begin{equation}\label{5.19}
	\beta = \frac{1-\tau}{2(L+1)}, \ \ \ \rho = \beta\ln\left(3+ \delta_\alpha^{-1}\right),
\end{equation} 
where $\delta_\alpha$ is so small that $\rho \geq r_3(N,D,E,m,\sigma)$. Then due to (\ref{5.18}), we have that
\begin{equation}\label{5.20}
	\begin{aligned}
		\|v_1 - v_2\|_{\mathbb{L}^{\infty}(D)} \leq 
		c_7(D,\sigma) \left(\beta\ln\left(3+ \delta_\alpha^{-1}\right)\right)^{d/m} 
		\left(3+ \delta_\alpha^{-1}\right)^{2\beta(L+1)} \delta_\alpha + \\
		+  c_8(N,D,E,m) \left(\beta\ln\left(3+ \delta_\alpha^{-1}\right)\right)^{-\frac{m-d}{m}} =\\
		= c_7(D,\sigma) \beta^{d/m} \left(1+ 3\delta_\alpha \right)^{1-\tau} \delta_\alpha^{\tau} 
		\left(\ln\left(3+ \delta_\alpha^{-1}\right)\right)^{d/m} + \\
		+ c_8(N,D,E,m) \beta^{-\frac{m-d}{m}}\left(\ln\left(3+ \delta_\alpha^{-1}\right)\right)^{-\frac{m-d}{m}},
	\end{aligned}
\end{equation}
where  $\tau, \beta$ and $\delta_\alpha$ are the same as in (\ref{5.19}). 

Using (\ref{5.20}), we obtain that
\begin{equation}\label{5.21}
	\|v_1 - v_2\|_{\mathbb{L}^{\infty}(D)} \leq c_{9}(N,D,E,m,\sigma)\left(\ln\left(3+\delta_\alpha ^{-1}\right)\right)^{-\frac{m-d}{m}}
\end{equation}
for $\delta_\alpha = \|\hat{M}_{\alpha, v_2} - \hat{M}_{\alpha, v_1}\|\leq \delta^{(0)}(N,D,E,m,\sigma)$, where
$\delta^{(0)}$ is a sufficiently small positive constant. Estimate (\ref{5.21}) in the general case (with modified $c_{9}$) follows 
from (\ref{5.21}) for $\delta_\alpha \leq \delta^{(0)}(N,D,E,m,\sigma)$ 
and the property that $\|v_j\|_{\mathbb{L}^{\infty}(D)} \leq c_{10}(D,m)N$.

	Thus, Theorem 2.1  is proved for $s = \frac{m-d}{m}$ and, since $\ln\left(3+ \delta_\alpha^{-1}\right) > 1$, for any  $0<s \leq \frac{m-d}{m}$. 
\section{Buckhgeim-type analogs of the Faddeev functions} 

In dimension $d=2$, we consider the functions $G_{z_0}, \psi_{z_0}, \tilde{\psi}_{z_0}$, $\delta h_{z_0}$ of \cite{NS2010}, 
going back to Buckhgeim's paper \cite{Buckhgeim2008} and being analogs of the Faddeev functions:
\begin{equation}\label{eq_bu_psi}
	\begin{aligned}
	\psi_{z_0}(z,\lambda) = e^{\lambda(z-z_0)^2} + \int\limits_{D} G_{z_0}(z,\zeta,\lambda)
	v(\zeta) \psi_{z_0}(\zeta,\lambda)\, d\mbox{Re}\zeta\, d\mbox{Im}\zeta,
	\\
	\widetilde{\psi}_{z_0}(z,\lambda) = e^{\bar{\lambda}(\bar{z}-\bar{z}_0)^2} + \int\limits_{D} \overline{G_{z_0}(z,\zeta,\lambda)}
	v(\zeta) \widetilde{\psi}_{z_0}(\zeta,\lambda)\, d\mbox{Re}\zeta\, d\mbox{Im}\zeta,
	\end{aligned}
\end{equation}
\begin{equation}
	\begin{aligned}\label{Bu_G}
		G_{z_0}(z,\zeta,\lambda) = \frac{1}{4\pi^2}\int\limits_{D} 
		\frac
		{e^{-\lambda(\eta-z_0)^2+\bar{\lambda}(\bar{\eta}-\bar{z_0})^2}d\mbox{Re}\eta\,d\mbox{Im}\eta}
		{(z-\eta)(\bar{\eta}-\bar{\zeta})} \, e^{\lambda(z-z_0)^2-\bar{\lambda}(\bar{\zeta}-\bar{z_0})^2},\\
		z = x_1 + ix_2,\  z_0\in D,\  \lambda \in \mathbb{C},
	\end{aligned}
\end{equation} 
where $\mathbb{R}^2$ is identified with $\mathbb{C}$ and $v$, $D$ satisfy (\ref{eq_c}), (\ref{eq_c2}) for $d=2$;
\begin{equation}\label{deltah}
	\delta h_{z_0}(\lambda) = \int\limits_{D} \widetilde{\psi}_{z_0,1}(z,-\lambda) 
	\left(v_2(z) - v_1(z)\right)
	\psi_{z_0,2}(z,\lambda) \, d\mbox{Re}z\, d\mbox{Im}z, \ \ \lambda \in \mathbb{C},
\end{equation}
where $v_1$, $v_2$ satisfy (\ref{eq_c2}) for $d=2$ and $\widetilde{\psi}_{z_0,1}$, $\psi_{z_0,2}$ denote 
$\widetilde{\psi}_{z_0}$, $\psi_{z_0}$ of (\ref{eq_bu_psi}) for $v=v_1$ and $v=v_2$, respectively.

We recall that (see \cite{NS2010}, \cite{NS2011}):
\begin{equation}\label{5.4b}
	\begin{aligned}
	4\frac{\partial^2}{\partial z \partial \bar{z}} G_{z_0}(z,\zeta,\lambda) = \delta(z - \zeta),\\
	4\frac{\partial^2}{\partial \zeta \partial \bar{\zeta}} G_{z_0}(z,\zeta,\lambda) = \delta(z - \zeta),
	\end{aligned}
\end{equation}
where $z,z_0, \zeta \in D$, $\lambda \in \mathbb{C}$ and $\delta$ is the Dirac delta function; formulas (\ref{eq_bu_psi}) 
at fixed $z_0$ and $\lambda$ are
considered as equations for $\psi_{z_0}$, $\widetilde{\psi}_{z_0}$ in $L^{\infty}(D)$; as a corollary of (\ref{eq_bu_psi}),
(\ref{Bu_G}), (\ref{5.4b}), the functions $\psi_{z_0}$, $\widetilde{\psi}_{z_0}$ satisfy (\ref{eq}) for $E=0$ and $d=2$;
$\delta h_{z_0}$ is similar to the right side of (\ref{4.15}). 

Let potentials $v,v_1,v_2 \in C^2(\bar{D})$ and 
\begin{equation}
	\begin{aligned}
	\|v\|_{C^2(\bar{D})}  \leq N, \ \ \|v_j\|_{C^2(\bar{D})}  \leq N, \ \ j=1,2,\\
	(v_1-v_2)|_{\partial D} = 0, \ \  \ \frac{\partial }{\partial \nu}(v_1-v_2)|_{\partial D} = 0,
	\end{aligned}
\end{equation}
then we have that:
\begin{equation}\label{5.5b}
	\begin{aligned}
		\psi_{z_0}(z,\lambda) = e^{\lambda(z-z_0)^2} \mu_{z_0}(z,\lambda), \ \ \
		\widetilde{\psi}_{z_0}(z,\lambda) = e^{\bar{\lambda}(\bar{z}-\bar{z}_0)^2} \widetilde{\mu}_{z_0}(z,\lambda),
	\end{aligned}
\end{equation} 
\begin{equation}\label{5.6b}
	\mu_{z_0}(z,\lambda) \rightarrow 1, \ \ \widetilde{\mu}_{z_0}(z,\lambda) \rightarrow 1 \ \ \text{ as } |\lambda| \rightarrow \infty
\end{equation}
and, for any $\sigma>1$,
\begin{subequations} \label{5.7b}
	\begin{equation}
		|\mu_{z_0}(z,\lambda)| + |\nabla\mu_{z_0}(z,\lambda)| \leq \sigma, 
	\end{equation}
	\begin{equation}
	|\widetilde{\mu}_{z_0}(z,\lambda)| + |\nabla\widetilde{\mu}_{z_0}(z,\lambda)|  \leq \sigma,
	\end{equation}
\end{subequations} 
where $\nabla = \left(\partial/\partial x_1, \partial/\partial x_2\right)$, $z= x_1 + ix_2$, $z_0 \in D$, $\lambda \in \mathbb{C}$, $|\lambda| \geq \rho_1(N,D,\sigma)$;
\begin{equation}\label{5.8b}
	\begin{aligned}
	v_2(z_0) - v_1(z_0) = \lim\limits_{\lambda\rightarrow \infty} \frac{2}{\pi}|\lambda| \delta h_{z_0}(\lambda)\\
		\text{ for any } z_0\in D,
	\end{aligned}
\end{equation}
\begin{equation}\label{5.9b}
\begin{aligned}
	\left| v_2(z_0) - v_1(z_0) -  \frac{2}{\pi}|\lambda| \delta h_{z_0}(\lambda) \right| \leq 
	\frac{c_{11}(N,D) \left(\ln(3|\lambda|)\right)^2}
	{|\lambda|^{3/4}}\\
		\text{ for } z_0\in D, \ |\lambda| \geq \rho_2(N,D). 
\end{aligned}
\end{equation}
Formulas (\ref{5.5b}) can be considered as definitions of $\mu_{z_0}$, $\widetilde{\mu}_{z_0}$. Formulas (\ref{5.6b}), (\ref{5.8b})
were given in \cite{NS2010}, \cite{NS2011} and go back to \cite{Buckhgeim2008}. Estimate (\ref{5.9b}) was obtained in \cite{NS2010}, \cite{S2011}. Estimates \eqref{5.7b} are proved in Section 8.
\section{Proof of Theorem \ref{Theorem_2.2}}
We suppose that 
 $\widetilde{\psi}_{z_0,1}(\cdot,-\lambda)$, $\psi_{z_0,2}(\cdot,\lambda)$, $\delta h_{z_0}(\lambda)$
are defined as in Section 6 but with $v_j-E$ in place of $v_j$, $j=1,2$. We use the identity
\begin{equation}
	\hat{M}_{\alpha,v}(E) = \hat{M}_{\alpha,v-E}(0).
\end{equation}
We also use the notation $N_E = N+E$. Then, using (\ref{5.9b}), we have that
\begin{equation}\label{7.1}
\begin{aligned}
	\left| v_2(z_0) - v_1(z_0) -  \frac{2}{\pi}|\lambda| \delta h_{z_0}(\lambda) \right| \leq 
	\frac{c_{11}(N_E,D) \left(\ln(3|\lambda|)\right)^2}
	{|\lambda|^{3/4}}\\
		\text{ for } z_0\in D, \ |\lambda| \geq \rho_2(N_E,D). 
\end{aligned}
\end{equation}
According to Theorem \ref{Th_ident} and (\ref{deltah}),  we get that
\begin{equation}\label{7.2}
			\begin{aligned}
	\delta h_{z_0}(\lambda)  =  \frac{1}{4\pi^2} \int\limits_{\partial D} [ \widetilde{\psi}_{z_0,1}(\cdot,-\lambda)]_\alpha  
	\left(\hat{M}_{\alpha,v_2}(E) - \hat{M}_{\alpha,v_1}(E)\right) 
	\left[\psi_{z_0,2}(\cdot,\lambda)\right]_\alpha   |dz|, \\ 
		\lambda \in \mathbb{C}.
	\end{aligned}
\end{equation}
Let 
\begin{equation}
	\begin{aligned}
	c_{12} = \frac{1}{4\pi^2}\int\limits_{\partial D} |dz|, \ \ \ L = \max\limits_{z\in \partial D}|z|,\\
	\delta_\alpha = \|\hat{M}_{\alpha, v_2}(E) - \hat{M}_{\alpha, v_1}(E)\|,
	\end{aligned}
\end{equation}
where $\|\hat{M}_{\alpha, v_2}(E) - \hat{M}_{\alpha, v_1}(E)\|$ is defined according to (\ref{norm_A}).

Using (\ref{7.2}), we get that
\begin{equation}\label{7.4}
	\begin{aligned}
	|\delta h_{z_0}(\lambda)| \leq  c_{12}
	\|[ \widetilde{\psi}_{z_0,1}(\cdot,-\lambda)]_\alpha  \|_{\mathbb{L}^{\infty}(\partial D)}\, \delta_\alpha \,
	\|\left[\psi_{z_0,2}(\cdot,\lambda)\right]_\alpha \|_{\mathbb{L}^{\infty}(\partial D)},\ \
		\lambda \in \mathbb{C}.
	\end{aligned}
\end{equation}
Using (\ref{[psi]_alpha}), (\ref{5.7b}), we find that:
\begin{equation}\label{7.5}
\begin{aligned}
	\|[\widetilde{\psi}_{z_0,1}(\cdot,-\lambda)]_\alpha\|_{\mathbb{L}^{\infty}(\partial D)} \leq
	 \sigma \,\exp\bigg(|\lambda|(4L^2+4L)\bigg), \\ 
	 \|[\psi_{z_0,2}(\cdot,\lambda)]_\alpha\|_{\mathbb{L}^{\infty}(\partial D)} \leq
	 \sigma \,\exp\bigg(|\lambda|(4L^2+4L)\bigg), 
	 \\
	\lambda \in \mathbb{C}, \ \  |\lambda|\geq\rho_1(N_E,D,\sigma).
\end{aligned}
\end{equation}
Here and bellow in this section the constant $\sigma$ is the same that in (\ref{5.7b}).

Combining (\ref{7.4}), (\ref{7.5}), we obtain that
\begin{equation}\label{7.6}
	\begin{aligned}
	|\delta h_{z_0}(\lambda)| \leq  c_{12} \sigma^2 \,\exp\bigg(|\lambda|(8L^2+8L)\bigg) \delta_\alpha, \\
		\lambda \in \mathbb{C}, \ \  |\lambda|\geq\rho_1(N_E,D,\sigma).
	\end{aligned}
\end{equation}
Using (\ref{7.1}) and (\ref{7.6}), we get that
\begin{equation}\label{7.7}
	\begin{aligned}
		\left| v_2(z_0) - v_1(z_0) \right| \leq  c_{12} \sigma^2 \,\exp\bigg(|\lambda|(8L^2+8L)\bigg) \delta_\alpha 
		+ 	\frac{c_{11}(N_E,D) \left(\ln(3|\lambda|)\right)^2}
	{|\lambda|^{3/4}}, \\
		z_0 \in D, \ \ \lambda \in \mathbb{C}, \ \  |\lambda|\geq\rho_3(N_E,D,\sigma) = \max\{\rho_1,\rho_2\}.
	\end{aligned}
\end{equation}

We fix some $\tau \in (0,1)$ and let
\begin{equation}\label{7.8}
	\beta = \frac{1-\tau}{8L^2+8L}, \ \ \ \lambda = \beta\ln\left(3+ \delta_\alpha^{-1}\right),
\end{equation} 
where $\delta_\alpha$ is so small that $|\lambda| \geq \rho_3(N_E,D,\sigma)$. Then due to (\ref{7.7}), we have that
\begin{equation}\label{7.9}
	\begin{aligned}
		\|v_1 - v_2\|_{\mathbb{L}^{\infty}(D)} &\leq 
		c_{12} \sigma^2 
		\left(3+ \delta_\alpha^{-1}\right)^{\beta(8L^2+8L)} \delta_\alpha + \\
	&+  c_{11}(N_E,D)  	\frac{
	\left(\ln\left(3\beta\ln\left(3+ \delta_\alpha^{-1}\right)\right)\right)^2}
	{\left(\beta\ln\left(3+ \delta_\alpha^{-1}\right)\right)^{\frac{3}{4}} }  =\\
		&= c_{12} \sigma^2 \left(1+ 3\delta_\alpha \right)^{1-\tau} \delta_\alpha^{\tau} 
		 + \\
		&+ c_{11}(N_E,D)  \beta^{-\frac{3}{4}}	\frac{
	\left(\ln\left(3\beta\ln\left(3+ \delta_\alpha^{-1}\right)\right)\right)^2}
	{\left(\ln\left(3+ \delta_\alpha^{-1}\right)\right)^{\frac{3}{4}} },
	\end{aligned}
\end{equation}
where  $\tau, \beta$ and $\delta_\alpha$ are the same as in (\ref{7.8}). 

Using (\ref{7.9}), we obtain that
\begin{equation}\label{7.10}
	\|v_1 - v_2\|_{\mathbb{L}^{\infty}(D)} \leq c_{13}(N_E,D,\sigma)\left(\ln\left(3+\delta_\alpha ^{-1}\right)\right)^{-\frac{3}{4}} 
	\left(\ln\left(3\ln\left(3+ \delta_\alpha^{-1}\right)\right)\right)^2
\end{equation}
for $\delta_\alpha = \|\hat{M}_{\alpha, v_2}(E) - \hat{M}_{\alpha, v_1}(E)\|\leq \delta^{(0)}(N_E,D,\sigma)$, where
$\delta^{(0)}$ is a sufficiently small positive constant. Estimate (\ref{5.21}) in the general case (with modified $c_{13}$) follows 
from (\ref{7.10}) for $\delta_\alpha \leq \delta^{(0)}(N_E,D,\sigma)$ 
and the property that $\|v_j\|_{\mathbb{L}^{\infty}(D)} \leq c_{14}(D)N$.

	Thus, Theorem 2.2  is proved for $s = \frac{3}{4}$ and, since $\ln\left(3+ \delta_\alpha^{-1}\right) > 1$, for any  $0<s \leq \frac{3}{4}$. 
	
\section{Proof of estimates (\ref{5.7b})}	

In this section we prove estimate (\ref{5.7b}a). Estimate  (\ref{5.7b}b) can be proved a completely similar way.
Let 
\begin{equation}
	\begin{aligned}
		C_{\bar{z}}^1 (\bar{D}) = \left\{u: u,\frac{\partial u}{\partial \bar{z}} \in C(\bar{D})\right\},\\
			\|u\|_{C_{\bar{z}}^1 (\bar{D})} = \max \left(\|u\|_{C(\bar{D})}, \|\frac{\partial u}{\partial \bar{z}}\|_{C(\bar{D})}\right).
	\end{aligned}
\end{equation}
Due to estimates of Section 3 of \cite{NS2010}, we have that, for any $\varepsilon_1>0$,
\begin{equation}\label{8.2}
	\begin{aligned}
		\mu_{z_0}(\cdot,\lambda) \in C_{\bar{z}}^1 (\bar{D}), \ \ \ \
		\|\mu_{z_0}(\cdot,\lambda)\|_{C_{\bar{z}}^1 (\bar{D})} \leq 1 + \varepsilon_1\ \
		\text{ for } \ |\lambda| \geq \rho_4(N,D,\varepsilon_1).
	\end{aligned}
\end{equation}
In view of (\ref{8.2}), to prove (\ref{5.7b}a) it remains to prove that, for any $\varepsilon_2>0$, 
\begin{equation}\label{8.3}
	\begin{aligned}
		\partial_{z}\mu_{z_0}(\cdot,\lambda) \in C(\bar{D}), \ \ \ \
		\|\partial_{z}\mu_{z_0}(\cdot,\lambda)\|_{C (\bar{D})} \leq \varepsilon_2\ \
		\text{ for } \ |\lambda| \geq \rho_5(N,D,\varepsilon_2),
	\end{aligned}
\end{equation}
where $\partial_{z}\mu_{z_0}(\cdot,\lambda)$ is considered as a function of $z\in \bar{D}$ and $\partial_{z} = \partial/\partial z$.

We have that (see Sections 2 and 5 of \cite{NS2010}):
\begin{equation}\label{8.4}
	\partial_{z}\mu_{z_0} = \frac{1}{4}\Pi\, \bar{\mbox{T}}_{z_0,\lambda} v \mu_{z_0},
\end{equation}
\begin{equation}
	\Pi u(z) = -\frac{1}{\pi} \int\limits_{D} \frac{u(\zeta)}{(\zeta-z)^2} d\mbox{Re}\zeta\, d\mbox{Im}\,\zeta,
\end{equation}
\begin{equation}
	\bar{\mbox{T}}_{z_0,\lambda} u(z) = -
	\frac{e^{-\lambda(z-z_0)^2+\bar{\lambda}(\bar{z}-\bar{z_0})^2}}{\pi}
	\int\limits_{D} \frac{e^{\lambda(\zeta-z_0)^2-\bar{\lambda}(\bar{\zeta}-\bar{z_0})^2}}{\bar{\zeta}-\bar{z}} u(\zeta)
	d\mbox{Re}\zeta\, d\mbox{Im}\,\zeta,
\end{equation}
where $u$ is a test function, $z\in \bar{D}$.

In view of (\ref{8.2}), (\ref{8.4}) and Theorem 1.33 of \cite{Vekua1962}, to prove (\ref{8.3}) it is sufficient to show  
that 
\begin{equation}\label{8.7}
	\|\bar{\mbox{T}}_{z_0,\lambda} u\|_{C_s(\bar{D})} \leq \frac{A(D,s)}{|\lambda|^{\delta(s)}}
	||u||_{C^1_{\bar{z}}(\bar{D})}, \ \ \ |\lambda|\geq 1, \ \ z_0 \in \bar{D},
\end{equation}
for some fixed $s\in (0,\frac{1}{2})$ and $\delta(s)>0$,
where $C_s(\bar{D})$ is the H\"older space, 
\begin{equation}
	\begin{aligned}
		&C_s(\bar{D}) = \left\{ u\in C(\bar{D}): \|u\|_{C_s(\bar{D})}<+\infty\right\},\\
		&\|u\|_{C_s(\bar{D})} = \max \left\{ \|u\|_{C(\bar{D})}, \|u\|'_{C_s(\bar{D})}\right\},\\
	&\|u\|'_{C_s(\bar{D})} = \sup_{z_1,z_2\in \bar{D}, 0<|z_1 - z_2|<1} \frac{|u(z_1) - u(z_2)|}{|z_1-z_2|^s}.
	\end{aligned}
\end{equation}
Due to estimate (5.6) of \cite{NS2010}, we have that
\begin{equation}\label{8.9}
	\|\bar{\mbox{T}}_{z_0,\lambda} u\|_{C(\bar{D})} \leq \frac{A_0(D)}{|\lambda|^{1/2}}
	||u||_{C^1_{\bar{z}}(\bar{D})}, \ \ \ |\lambda|\geq 1, \ \ z_0 \in \bar{D}.
\end{equation}
Therefore, to prove \eqref{8.7} it remains to prove that
\begin{equation}\label{8.10}
	\|\bar{\mbox{T}}_{z_0,\lambda} u\|'_{C_s(\bar{D})} \leq \frac{A_1(D,s)}{|\lambda|^{\delta(s)}}
	||u||_{C^1_{\bar{z}}(\bar{D})}, \ \ \ |\lambda|\geq 1, \ \ z_0 \in \bar{D}, 
\end{equation} 
for some fixed $s\in (0,\frac{1}{2})$ and $\delta(s)>0$.

We will use that
\begin{equation}\label{8.11}
	\|u_1u_2\|'_{C_s(\bar{D})}\leq \|u_1\|'_{C_s(\bar{D})} \|u_2\|_{C(\bar{D})} +
	\|u_1\|_{C(\bar{D})} \|u_2\|'_{C_s(\bar{D})}, \ \ \ 0<s<1.
\end{equation}
One can see that
\begin{equation}\label{8.12}
	\bar{\mbox{T}}_{z_0,\lambda} = \mbox{F}_{z_0,-\lambda}
	\bar{\mbox{T}} {\mbox{F}}_{z_0,\lambda},
\end{equation}
where $\bar{\mbox{T}} = \bar{\mbox{T}}_{z_0,0}$ and ${\mbox{F}}_{z_0,\lambda}$ is the 
multiplication operator by the function 
\begin{equation}\label{8.13}
	F(z,z_0,\lambda) = {e^{\lambda(z-z_0)^2-\bar{\lambda}(\bar{z}-\bar{z_0})^2}}.
\end{equation}

One can see also that
\begin{equation}\label{8.14}
	\begin{aligned}
		&\|F(\cdot,z_0,-\lambda)\|_{C(\bar{D})} = 1,\\
		&\|F(\cdot,z_0,-\lambda)\|'_{C_s(\bar{D})} \leq A_2(D,s) |\lambda|^s, \ \ \ |\lambda|\geq 1, \ \ z_0 \in \bar{D}.
	\end{aligned}
\end{equation}
In view of \eqref{8.9}, \eqref{8.11} - \eqref{8.14}, to prove \eqref{8.10} it remains to prove that
\begin{equation}\label{8.15}
	\|\bar{\mbox{T}}{\mbox{F}}_{z_0,\lambda} u\|'_{C_s(\bar{D})} \leq \frac{A_3(D,s)}{|\lambda|^{\delta_1(s)}}
	||u||_{C^1_{\bar{z}}(\bar{D})}, \ \ \ |\lambda|\geq 1, \ \ z_0 \in \bar{D}, 
\end{equation} 
for some fixed $s\in (0,\frac{1}{2})$ and $\delta_1(s)>0$.

We have that
\begin{equation}
	\begin{aligned}
	\pi\bar{\mbox{T}}{\mbox{F}}_{z_0,\lambda} u(z_1) - \pi\bar{\mbox{T}}{\mbox{F}}_{z_0,\lambda} u(z_2) 
	=
	\int\limits_{D} \frac{F(\zeta,z_0,\lambda) u(\zeta) (\bar{z}_2-\bar{z}_1)}
	{(\bar{\zeta}-\bar{z}_1)(\bar{\zeta}-\bar{z}_2)} 
	d\mbox{Re}\zeta\, d\mbox{Im}\,\zeta = \\
	= I_{z_0,\lambda,\varepsilon}(z_1,z_2) + J_{z_0,\lambda,\varepsilon}(z_1,z_2),
	\end{aligned}
\end{equation}
where 
\begin{equation}
	I_{z_0,\lambda,\varepsilon}(z_1,z_2) =
	\int\limits_{D\setminus D_{z_0,z_1,z_2,\varepsilon}} \frac{F(\zeta,z_0,\lambda) u(\zeta) (\bar{z}_2-\bar{z}_1)}
	{(\bar{\zeta}-\bar{z}_1)(\bar{\zeta}-\bar{z}_2)} 
	d\mbox{Re}\zeta\, d\mbox{Im}\,\zeta,
\end{equation}
\begin{equation}
	J_{z_0,\lambda,\varepsilon}(z_1,z_2) =
	\int\limits_{D_{z_0,z_1,z_2,\varepsilon}} \frac{F(\zeta,z_0,\lambda) u(\zeta) (\bar{z}_2-\bar{z}_1)}
	{(\bar{\zeta}-\bar{z}_1)(\bar{\zeta}-\bar{z}_2)} 
	d\mbox{Re}\zeta\, d\mbox{Im}\,\zeta,
\end{equation}
where $B_{z,\varepsilon} = \{\zeta \in \mathbb{C}: |\zeta-z|<\varepsilon\}$,
$D_{z_0,z_1,z_2,\varepsilon} = D \setminus \left(\bigcup\limits_{j=0}\limits^{2} B_{z_j,\varepsilon}\right)$.

We will use the following inequalities:
\begin{align}
&\left|\frac{z_2-z_1}{(\zeta-z_1)(\zeta-z_2)}\right| \leq n_1 |z_2-z_1|^{s} \sum\limits_{j=1}\limits^{2} \frac{1}{|\zeta-z_j|^{1+s}},
\\
&\left|\frac{z_2-z_1}{(\zeta-z_1)(\zeta-z_2)(\zeta-z_0)}\right| \leq 
n_2 |z_2-z_1|^{s} \sum\limits_{j=0}\limits^{2} \frac{1}{|\zeta-z_j|^{2+s}},\\
&\left|\frac{\partial}{\partial \zeta}\left( \frac{z_2-z_1}{(\zeta-z_1)(\zeta-z_2)(\zeta-z_0)}\right)\right| \leq 
n_3 |z_2-z_1|^{s} \sum\limits_{j=0}\limits^{2} \frac{1}{|\zeta-z_j|^{3+s}},
\end{align}
where $s\in (0,1)$,  $n_1,n_2,n_3>0$,  $z_0,z_1,z_2,\zeta \in \mathbb{C}$ and $\zeta \neq z_i$ for $j=0,1,2$.

Using (8.17), (8.19), we obtain that
\begin{equation}
	I_{z_0,\lambda,\varepsilon}(z_1,z_2) \leq n_4(s) \varepsilon^{1-s}|z_2-z_1|^s,
\end{equation}
where $n_4(s)>0$, $z_0,z_1,z_2,\zeta \in \mathbb{C}$  and $\varepsilon \in (0,1)$. Further, we have that
\begin{equation}
\begin{aligned}
	J_{z_0,\lambda,\varepsilon}(z_1,z_2) = -\frac{1}{2\bar{\lambda}}
	\int\limits_{D_{z_0,z_1,z_2,\varepsilon}} \frac{\partial F(\zeta,z_0,\lambda)}{\partial \bar{\zeta}} 
	\frac{ u(\zeta) (\bar{z}_2-\bar{z}_1)}
	{(\bar{\zeta}-\bar{z}_1)(\bar{\zeta}-\bar{z}_2)(\bar{\zeta}-\bar{z}_0)} 
	d\mbox{Re}\zeta\, d\mbox{Im}\,\zeta=\\
	= J^1_{z_0,\lambda,\varepsilon}(z_1,z_2) + J^2_{z_0,\lambda,\varepsilon}(z_1,z_2),
\end{aligned}
\end{equation}
where 
\begin{equation}
\begin{aligned}
	&J^1_{z_0,\lambda,\varepsilon}(z_1,z_2) = -\frac{1}{4i\bar{\lambda}}
	\int\limits_{\partial D_{z_0,z_1,z_2,\varepsilon}} \frac{F(\zeta,z_0,\lambda) u(\zeta) (\bar{z}_2-\bar{z}_1)}
	{(\bar{\zeta}-\bar{z}_1)(\bar{\zeta}-\bar{z}_2)(\bar{\zeta}-\bar{z}_0)} 
	d\zeta,\\
	&J_{z_0,\lambda,\varepsilon}^2(z_1,z_2) =\frac{1}{2\bar{\lambda}}
	\int\limits_{D_{z_0,z_1,z_2,\varepsilon}} F(\zeta,z_0,\lambda) \frac{\partial}{\partial \bar{\zeta}} \left(
	\frac{ u(\zeta) (\bar{z}_2-\bar{z}_1)}
	{(\bar{\zeta}-\bar{z}_1)(\bar{\zeta}-\bar{z}_2)(\bar{\zeta}-\bar{z}_0)} \right)
	d\mbox{Re}\zeta\, d\mbox{Im}\,\zeta,
\end{aligned}
\end{equation}
Using (8.20), (8.21), (8.24), we obtain that
\begin{equation}
	\begin{aligned}
	J^1_{z_0,\lambda,\varepsilon}(z_1,z_2) \leq |\lambda|^{-1} &n_5(D,s) \varepsilon^{-1-s}|z_2-z_1|^s \|u\|_{C(\bar{D})},\\
	J^2_{z_0,\lambda,\varepsilon}(z_1,z_2) \leq |\lambda|^{-1} &n_6(D,s) \varepsilon^{-1-s}|z_2-z_1|^s \|u\|_{C(\bar{D})}+\\
	&+|\lambda|^{-1} n_7(D,s) \varepsilon^{-s}|z_2-z_1|^s \left\|\frac{\partial u}{\partial \bar{z}}\right\|_{C(\bar{D})},
	\end{aligned}
\end{equation}
where $z_0,z_1,z_2,\lambda \in \mathbb{C}$, $|\lambda|\geq 1$, $\varepsilon \in (0,1)$.

Using (8.16), (8.22), (8.23), (8.25) and putting $\varepsilon = |\lambda|^{-1/2}$ into (8.22), (8.25),
we obtain (8.15) with $\delta_1(s)=(1-s)/2$.

\section{Proof of symmetry (\ref{3.10})}

	Let $D'$ be an open bounded domain in $\mathbb{R}^d$ such that 
	\begin{itemize}
	 \item $D \subset D'$, 
	 \item $D'$ satisfies (1.2),
	 \item $E$ is not a Dirichlet eigenvalue for the operator $-\Delta + v$ in $D'$. 
	\end{itemize}
	Here and bellow in this section
	we assume that $v\equiv 0$ on $D'\setminus D$. Let $R(x,y,E)$ denote the Green function for the
	operator $-\Delta+v-E$ in $D'$ with the Dirichlet boundary condition.  
	 We recall that
	\begin{equation}\label{Symm_R}
		R(x,y,E) = R(y,x,E), \ \ \ x,y\in D'.
	\end{equation}
	
	Using (\ref{G_F}), (\ref{Symm_R}), we find that for  $x,y\in D$
	\begin{equation}\label{R-R}
		\begin{aligned}
			\int\limits_{\partial D} &\left( R(x,\xi,E)  \frac{\partial R}{\partial \nu_\xi} (y,\xi,E)
		- R(y,\xi,E)  \frac{\partial R}{\partial \nu_\xi} (x,\xi,E) \right)  d\xi = \\
		&=\int\limits_{ D} \Big( R(x,\xi,E)  \left(\Delta_\xi - v + E \right) R(y,\xi,E) -  
		R(y,\xi,E)  \left(\Delta_\xi - v + E \right) R(x,\xi,E) \Big)  d\xi =\\
		&\ \ \ \ \ \ \text{ } \ \ \ \  \text{ } \ \ \ \  \text{ } \ \ \ \  \text{ } \ \ \ \
		 \text{ } \ \ \ \ \text{ } \ \ \ \ \text{ } \ \ \ \  \text{ } \ \ \ \  \text{ } \ \ \ \  \text{ } \ \ \ \
		= - R(x,y,E) + R(y,x,E) = 0.
		\end{aligned}
	\end{equation}
	Note that $W = G_\alpha + R(E)$ is the solution of the equation 
	\begin{equation}\label{eq_W}
		(-\Delta_x+v - E) W(x,y) = 0, \ \ x,y\in D
	\end{equation} 
	with the boundary condition
	\begin{equation}\label{bc_W}
		\begin{aligned}
		&\left(\cos \alpha \,W(x,y) - \sin\alpha\, \frac{\partial W}{\partial \nu_x}(x,y)\right) \Big|_{x\in \partial D} = \\
		&= \left(\cos \alpha \,R(x,y,E) - \sin\alpha\, \frac{\partial R}{\partial \nu_x}(x,y,E)\right) \Big|_{x\in \partial D}, \ \ y\in D.
		\end{aligned}
	\end{equation}
	Using (\ref{G_F}) and (\ref{eq_W}), we find that for $x,y \in D$
	\begin{equation}\label{W-W}
		\begin{aligned}
			\int\limits_{\partial D} &\left( W(\xi,x)  \frac{\partial W}{\partial \nu_\xi} (\xi,y)
		- W(\xi,y) \frac{\partial W}{\partial \nu_\xi} (\xi,x) \right)  d\xi = \\
		&=\int\limits_{ D} \Big( W(\xi,x)  \left(\Delta_\xi - v + E \right) W(\xi,y) -  
		W(\xi,y)  \left(\Delta_\xi - v + E \right) W(\xi,x) \Big)  d\xi = 0
		\end{aligned}
	\end{equation}
	Note that
	\begin{equation}\label{eq_WR}
		W(x,y) = -\int\limits_{ D}  W(\xi,y)  \left(\Delta_\xi - v + E \right) R(\xi,x,E) d \xi, \ \ \ x,y\in D.
	\end{equation}
	Combining (\ref{G_F}), (\ref{eq_W}) and (\ref{eq_WR}), we obtain that
	\begin{equation}\label{part_WR}
		\begin{aligned}
		W(x,y) = 	-\int\limits_{\partial D} \left( W(\xi,y)  \frac{\partial R}{\partial \nu_\xi} (\xi,x,E)
		- R(\xi,x,E) \frac{\partial W}{\partial \nu_\xi} (\xi,y) \right)  d\xi, \\ x,y\in D.
		\end{aligned}
	\end{equation}
	Using (\ref{bc_W}) and (\ref{part_WR}), we get that
	\begin{equation}\label{sin_W}
		\begin{aligned}
			\sin\alpha\, &W(x,y) = \\ =\int\limits_{\partial D} &W(\xi,y) 
			\left(\cos\alpha \, W(\xi,x) - \sin\alpha \, \frac{\partial W}{\partial \nu_\xi} (\xi,x)  - \cos\alpha\, R(\xi,x,E)\right) d\xi 
			\ -\\ &-
			\int\limits_{\partial D} R(\xi,x,E) 
			\left(\cos\alpha \, R(\xi,y,E) - \sin\alpha \, \frac{\partial R}{\partial \nu_\xi} (\xi,x,E)  - \cos\alpha\, W(\xi,y)\right) d\xi,\\ 
			 &\ \ \ \ \ \ \text{ } \ \ \ \  \text{ } \ \ \ \  \text{ } \ \ \ \  \text{ } \ \ \ \
		 \text{ } \ \ \ \ \text{ } \ \ \ \ \text{ } \ \ \ \  \text{ } \ \ \ \  \text{ } \ \ \ \  \text{ } \ \ \ \  \text{ } \ \ \ \
		  \text{ } \ \ \ \ \text{ } \ \ \ \ \text{ } \ \ \ \ \text{ } \ \ \ \
			 x,y\in D.
		\end{aligned}
	\end{equation}  
	Combining similar to (\ref{sin_W}) formula for $\sin\alpha\, W(y,x)$, (\ref{R-R}) and (\ref{W-W}), we obtain that
	\begin{equation}
		\sin\alpha\, W(x,y) - \sin\alpha\, W(y,x) = 0, \ \ \ x,y\in D.
	\end{equation} 
	In the case of $\sin \alpha = 0$, combining (\ref{bc_W}) and (\ref{part_WR}), we get that
	\begin{equation}
		\begin{aligned}
			W(x,y) = 	\int\limits_{\partial D} \left( - R(\xi,y,E)  \frac{\partial R}{\partial \nu_\xi} (\xi,x,E)
		+ W(\xi,x) \frac{\partial W}{\partial \nu_\xi} (\xi,y) \right)  d\xi, \\ x,y\in D.
		\end{aligned}
	\end{equation}
	Hence, one can get that for any $\alpha$
	\begin{equation}\label{Symm_W}
		W(x,y) = W(y,x),  \ \ \ x,y\in D.
	\end{equation}
	Combining (\ref{Symm_R}) and (\ref{Symm_W}), we obtain (\ref{3.10}).

We note that symmetry (\ref{3.10}) for  $v \equiv 0$, $E=0$, $d\geq 3$ was proved early, for example,  in \cite{Lanzani2004}. 

\section*{Acknowledgements}
The second author was partially supported by the Russian Federation Goverment grant No. 2010-220-01-077.

\noindent
{ {\bf M.I. Isaev}\\
Centre de Math\'ematiques Appliqu\'ees, Ecole Polytechnique,

91128 Palaiseau, France\\
Moscow Institute of Physics and Technology,

141700 Dolgoprudny, Russia\\
e-mail: \tt{isaev.m.i@gmail.com}}\\

\noindent
{ {\bf R.G. Novikov}\\
Centre de Math\'ematiques Appliqu\'ees, Ecole Polytechnique,

91128 Palaiseau, France\\
Institute of Earthquake Prediction Theory and Mathematical Geophysics RAS,

117997 Moscow, Russia\\ 
e-mail: \tt{novikov@cmap.polytechnique.fr}}

\end{document}